\documentclass[a4paper,10pt]{article}
\usepackage{stmaryrd}
\usepackage{amsfonts}
\usepackage{bbm}
\usepackage{amscd}
\usepackage{mathrsfs}
\usepackage{latexsym,amssymb,amsmath,amscd,amscd,amsthm,amsxtra}
\usepackage[dvips]{graphicx}
\usepackage[utf8]{inputenc}
\usepackage[T1]{fontenc}
\usepackage{lmodern}
\usepackage{amssymb}
\usepackage[all]{xy}
\usepackage{nicefrac,mathtools,enumitem}
\usepackage{microtype}
\usepackage{xcolor}
\newtheorem{thm}{Theorem}[section]
\newtheorem{lem}[thm]{Lemma}
\newtheorem{cor}[thm]{Corollary}
\newtheorem{pro}[thm]{Proposition}
\newtheorem{ex}[thm]{Example}
\newtheorem{rmk}[thm]{Remark}
\newtheorem{defi}[thm]{Definition}

\newcommand {\emptycomment}[1]{}

\newcommand {\nrn}[1]{   [    #1   ]^{\rm{3Lie}}}

\newcommand{\be }{\begin{equation}}
\newcommand{\ee }{\end{equation}}

\newcommand{\pf}{\noindent{\bf Proof.}\ }

\newcommand{\g}{\mathfrak g}

\newcommand{\frkg}{\mathfrak g}

\newcommand{\frkX}{\mathfrak X}
\newcommand{\frkY}{\mathfrak Y}

\def\qed{\hfill ~\vrule height6pt width6pt depth0pt}

\newcommand{\half}{\frac{1}{2}}

\newcommand{\br}[1]{   [ \cdot,    \cdot  ]   }

\newcommand{\dM}{\mathrm{d}}

\newcommand{\Hom}{\mathrm{Hom}}

\newcommand{\gl}{\mathfrak {gl}}

\newcommand{\End}{\mathrm{End}}
\newcommand{\ad}{\mathrm{ad}}

\newcommand{\K}{\mathbb{K}}

\begin{document}
\title{
{A new approach to representations of $3$-Lie algebras\\ and abelian extensions
\thanks
 {
This research is supported by NSFC (11471139) and NSF of Jilin Province (20140520054JH).
 }
} }
\author{Jiefeng Liu$^1$, Abdenacer Makhlouf$^2$ and  Yunhe Sheng$^3$
\\
\\
$^1$ Department of Mathematics, Xinyang Normal University,\\ \vspace{2mm} Xinyang 464000, Henan, China \\
$^2$University of Haute Alsace, Laboratoire de Math\'ematiques,\\ \vspace{2mm}Informatique et Applications, Mulhouse, France\\\vspace{2mm}
$^3$Department of Mathematics, Jilin University,
 Changchun 130012,  China
\\
email:
liujf12@126.com, abdenacer.makhlouf@uha.fr, shengyh@jlu.edu.cn
 }
\date{}
\footnotetext{{\it{Keyword}:  $3$-Lie algebra, representation, generalized representation, cohomology, abelian extension, Maurer-Cartan element  }} \footnotetext{{\it{MSC}}: 17B10, 17B56, 17A42.}
\maketitle

\begin{abstract}
In this paper, we introduce the notion of  generalized representation of a $3$-Lie algebra, by which we obtain a generalized semidirect product $3$-Lie algebra.  Moreover, we develop the corresponding cohomology theory.  Various examples of generalized representations of 3-Lie algebras and computation of 2-cocycles of the new cohomology are provided. Also, we show that a split abelian extension of a  3-Lie algebra is isomorphic to a generalized semidirect product $3$-Lie algebra. Furthermore, we describe general abelian extensions of 3-Lie algebras using Maurer-Cartan elements.
\end{abstract}


\section{Introduction}
The notion of a Filippov algebra, or an $n$-Lie algebra was introduced in  \cite{Filippov}. Ternary Lie algebras are related  to Nambu mechanics (\cite{N}), generalizing Hamiltonian mechanics by using more than one hamiltonian. The algebraic formulation of this theory is due to Takhtajan (\cite{T}), see also \cite{Gautheron}. Moreover, 3-Lie algebras appeared in String Theory. Fuzzy sphere (noncommutative space) arises naturally in the description of D1-branes ending
on D3-branes in Type IIB superstring theory and
 the effective dynamics of this system is described
by the Nahm equations.  In order to find an appropriate description of the lift of this configuration to
M-theory, one can study supergravity solutions describing M2-branes ending on M5-branes. In \cite{Basu},
 Basu and Harvey suggested to replace the Lie algebra appearing in the Nahm equation by a 3-Lie algebra for the  lifted Nahm equations.
Furthermore, in the context of Bagger-Lambert-Gustavsson model of multiple
M2-branes, Bagger-Lambert  managed to construct, using a ternary bracket,  an $N=2$
 supersymmetric version of the worldvolume theory of the M-theory membrane, see \cite{BL0}. An extensive literatures are related to this pioneering work, see \cite{BL3,BL2,HHM,P}.

In mathematics, the first appearance of ternary operation went back to 19th century, where Cayley considered cubic matrices. 
 In these last decades, several algebraic aspects of 3-Lie algebras, or more generally, of
$n$-Lie algebras   were studied. See   \cite{AMS11, Realization,Baiclassification} for the construction, realization and classifications of 3-Lie algebras and $n$-Lie algebras.  In particular, the representation theory of $n$-Lie algebras was first introduced by Kasymov in \cite{Kasymov}. The adjoint representation  is defined by the ternary bracket in which two elements are fixed.  Through fundamental objects one may also represent a 3-Lie algebra and more generally an $n$-Lie algebra by  a Leibniz algebra (\cite{DT}). Following this approach,  deformations of 3-Lie algebras and $n$-Lie algebras are studied in   \cite{deformation,Tcohomology}, see  \cite{Makhlouf} for a review. In \cite{NR bracket of n-Lie}, the author defined a graded Lie algebra structure on the cochain complex of an $n$-Leibniz algebra and describe an $n$-Leibniz structure as a canonical structure.  See  the review article \cite{review} for more details.

In this paper, we provide a new approach to representation theory of 3-Lie algebras. We define  generalized representations of a 3-Lie algebra $\g$ on a vector space $V$ via canonical structures in the differential graded Lie algebra associated to $\g\oplus V$ given in \cite{DT}. A generalized representation  leads also to a new 3-Lie algebra, which we call a generalized semidirect product. We also develop the corresponding cohomology theory. As applications, we study abelian extensions of 3-Lie algebras and provide several examples. We show that a split abelian extension is isomorphic to a generalized semidirect product. Note that even split abelian extensions cannot be studied via the usual representations (in the sense of Kasymov (\cite{Kasymov})). This justifies the usage of our approach. Furthermore, we use Maurer-Cartan elements to describe non-split abelian extensions.

The paper is organized as follows. In Section 2, we give a review of representations of 3-Lie algebras, the differential graded Lie algebras associated to 3-Lie algebras and canonical structures. In Section 3, we first explain  the usual representations of a 3-Lie algebra using canonical structures,  then we introduce the notion of  generalized representation of a 3-Lie algebra. We show that a generalized representation of a 3-Lie algebra  gives rise to a   generalized semidirect product 3-Lie algebra. Various examples of generalized representations are given. In Section 4, we provide the corresponding cohomology theory associated to generalized representations. Moreover, we study  in detail $1$-cocycles, $2$-cocycles and provide examples. 
Section 5 deals with abelian extensions of 3-Lie algebras. First,  we show that a split abelian extension is isomorphic to a generalized semidirect product. Unlike the case of Lie algebras, one cannot use cocycles to describe general abelian extensions of 3-Lie algebras. Alternatively, we can use Maurer-Cartan elements to describe them.

\section{$3$-Lie algebras and their representations}
In this paper, we work over an algebraically closed field $\K$ of characteristic 0 and all the vector spaces are over $\K$.
\begin{defi}{\rm(\cite{Filippov})}
  A {\bf 3-Lie algebra}  is a vector space $\g$ together with a skew-symmetric linear map  $[\cdot,\cdot,\cdot]:
\otimes^3 \g\rightarrow \g$ such that the following {\bf Fundamental Identity (FI)} holds:
\begin{eqnarray}
\nonumber &&F_{x_1,x_2,x_3,x_4,x_5}\\
\nonumber&\triangleq&[x_1,x_2,[x_3,x_4,x_5]]-[[x_1,x_2,x_3],x_4,x_5]-[x_3,[x_1,x_2,x_4],x_5]-[x_3,x_4,[x_1,x_2,x_5]]\\
\label{eq:de1}&=&0.
\end{eqnarray}
\end{defi}

 Elements in $\wedge^2\g$ are called {\bf fundamental objects} of the $3$-Lie algebra $(\g,[\cdot,\cdot,\cdot])$. There is a bilinear operation $[\cdot,\cdot]_{\rm F}$ on $  \wedge^{2}\g$, which is given by
\begin{equation}\label{eq:bracketfunda}
~[\frkX,\frkY]_{\rm F}=[x_1,x_2,y_1]\wedge y_2+y_1\wedge[x_1,x_2,y_2],\quad \forall \frkX=x_1\wedge x_2, ~\frkY=y_1\wedge y_2.
\end{equation}
It is well-known that $(\wedge^2\g,[\cdot,\cdot]_{\rm F})$ is a Leibniz algebra (\cite{DT}), which plays  an important role in the theory of 3-Lie algebras.

\begin{defi}{\rm (\cite{Kasymov})}\label{defi:usualrep}
 A {\bf representation} $\rho$ of a $3$-Lie algebra $\frkg$ on a vector space  $V$ is a linear map $\rho:\wedge^2\frkg\longrightarrow \End( V),$ such that
\begin{eqnarray*}
  \rho(x_1,x_2)\rho(x_3,x_4)&=&\rho([x_1,x_2,x_3],x_4)+\rho(x_3,[x_1,x_2,x_4])+\rho(x_3,x_4)\rho(x_1,x_2),\\
\rho(x_1,[x_2,x_3,x_4])&=&\rho(x_3,x_4)\rho(x_1,x_2)-\rho(x_2,x_4)\rho(x_1,x_3)+\rho(x_2,x_3)\rho(x_1,x_4).
\end{eqnarray*}
\end{defi}

It is straightforward to obtain
\begin{lem}\label{lem:semidirectp}
Let $\g$ be a $3$-Lie algebra, $V$  a vector space and $\rho:
\wedge^2\g\rightarrow \gl(V)$  a skew-symmetric linear
map. Then $(V;\rho)$ is a representation of $\g$ if and only if there
is a $3$-Lie algebra structure $($called the semidirect product$)$
on the direct sum of vector spaces  $\g\oplus V$, defined by
\begin{equation}\label{eq:sum}
[x_1+v_1,x_2+v_2,x_3+v_3]_{\rho}=[x_1,x_2,x_3]+\rho(x_1,x_2)v_3+\rho(x_3,x_1)v_2+\rho(x_2,x_3)v_1,
\end{equation}
for $x_i\in \g, v_i\in V, 1\leq i\leq 3$. We denote this semidirect product $3$-Lie algebra by $\g\ltimes_\rho V.$
\end{lem}

A $p$-cochain on  $\frkg$ with  coefficients in a representation   $(V;\rho)$ is a linear map $$\alpha:\wedge^2\frkg\otimes\stackrel{(p-1)}{\cdots}\otimes\wedge^2\frkg\wedge\frkg\longrightarrow V.$$
Denote the space of $p$-cochains by $C^{p-1}(\g,V).$ The coboundary operator $\delta_\rho:C^{p-1}(\g,V)\longrightarrow C^{p}(\g,V)$  is given by
\begin{eqnarray}
\nonumber&&(\delta_\rho\alpha)(\frkX_1,\cdots ,\frkX_p,z)\\
\nonumber&=& \sum_{1\leq j<k}(-1)^j\alpha(\frkX_1,\cdots ,\hat{\frkX}_j,\cdots ,\frkX_{k-1},[\frkX_j,\frkX_k]_{\rm F},\frkX_{k+1},\cdots ,\frkX_{p},z)\\
\nonumber&&+\sum_{j=1}^p(-1)^j\alpha(\frkX_1,\cdots ,\hat{\frkX}_j,\cdots ,\frkX_{p},[\frkX_j,z])\\
\nonumber&&+\sum_{j=1}^p(-1)^{j+1}\rho(\frkX_j)\alpha(\frkX_1,\cdots ,\hat{\frkX}_j,\cdots ,\frkX_{p},z)\\
\label{eq:drho}&&+(-1)^{p+1}\Big(\rho(y_{p},z)\alpha(\frkX_1,\cdots ,\frkX_{p-1},x_{p} ) +\rho(z,x_{p})\alpha(\frkX_1,\cdots ,\frkX_{p-1},y_{p} ) \Big),
\end{eqnarray}
for all   $\frkX_i=(x_i,y_i)\in\wedge^2\frkg$ and $z\in\frkg.$ An element $\alpha\in {C}^{p-1}(\g,V)$ is called a $p$-cocycle if $\delta_\rho\alpha=0$; It is called a $p$-coboundary if there exists some $\beta\in {C}^{p-2}(\g,V)$ such that $\alpha=\delta_\rho\beta$. Denote by $Z^p(\g;V)$ and $B^p(\g;V)$ the set of $p$-cocycles and the set of $p$-coboundaries respectively. Then  the $p$-th cohomology group is \begin{equation}\label{eq:cohomology}H^p(\g;V)=Z^p(\g;V)/B^p(\g;V).\end{equation}

In \cite{NR bracket of n-Lie}, the author constructed a graded Lie algebra structure by which one can describe an $n$-Leibniz algebra structure as a canonical structure. Here, we give the precise formulas for the 3-Lie algebra case.

  Set $L_p=C^p(\g,\g)=\Hom(\wedge^2\frkg\otimes\stackrel{(p)}{\cdots }\otimes\wedge^2\frkg\wedge\frkg,\g)$ and $L=\oplus_{p\geq0} L_p$. Let $\alpha\in C^p(\g,\g),\beta\in C^q(\g,\g),\quad p,q\geq 0$. Let $\frkX_i=x_i\wedge y_i\in \wedge^2\g$ for $i=1,2,\cdots ,p+q$ and $x\in\g$. For each subset 
 $J=\{j_1,\cdots,j_{q+1}\}_{j_1<\cdots<j_{q+1}}\subset N\triangleq \{1,2,\cdots,p+q+1\}$,
  let $I=\{i_1,\cdots,i_p\}_{i_1<\cdots<i_p}=N/J$. Then we have

  \begin{thm}{\rm (\cite{NR bracket of n-Lie})}\label{thm:gradelie}
  The graded vector space $L$ equipped with the graded commutator bracket
\begin{eqnarray}
\nrn{\alpha,\beta}=(-1)^{pq} \alpha\circ \beta-\beta\circ \alpha,
\end{eqnarray}
is a graded Lie algebra and $\alpha\circ \beta\in L^{p+q}$ is defined by
\begin{eqnarray*}
&&\alpha\circ \beta(\frkX_1,\cdots ,\frkX_{p+q},x)=\sum_{J,j_{q+1}<p+q+1}(-1)^{(J,I)}(-1)^k\\
&&\Big(\alpha(\frkX_{i_1},\cdots ,\frkX_{i_k},\beta(\frkX_{j_1},\cdots ,\frkX_{j_q},x_{j_{q+1}})\wedge y_{j_{q+1}},\frkX_{i_{k+1}},\cdots ,X_{i_{p-1}},x)\\
&&+\alpha(\frkX_{i_1},\cdots ,\frkX_{i_k},x_{j_{q+1}}\wedge \beta(\frkX_{j_1},\cdots ,\frkX_{j_q},y_{j_{q+1}}) ,\frkX_{i_{k+1}},\cdots ,X_{i_{p-1}},x)\Big)\\
&&+\sum_{J,j_{q+1}=p+q+1}(-1)^{(J,I)}(-1)^p\alpha(\frkX_{i_1},\cdots ,\frkX_{i_p},\beta(\frkX_{j_1},\cdots ,\frkX_{j_q},x)),
\end{eqnarray*}
where $k$ is uniquely determined by the condition  $i_k<j_{q+1}<i_{k+1}$ and if $j_{q+1}<i_1$, i.e. $j_{q+1}=q+1, i_1=q+2$ then $k=0$; if $j_{q+1}>i_{p-1}$, i.e. $j_{q+1}=p+q$ then $k=p-1$.
\end{thm}

We can use the graded Lie algebra structure $(L,\nrn{\cdot,\cdot})$ to describe 3-Lie algebra structures as well as coboundary operators.

\begin{lem}\label{lem:cs} The map
 $\pi:\wedge^3\g\longrightarrow\g$ defines a $3$-Lie bracket if and only if $\nrn{\pi,\pi}=0$, i.e. $\pi$ is a canonical structure.
\end{lem}
Let $\g$ be a $3$-Lie algebra. Given $x_1,x_2\in\g$, define $\ad:\wedge^2\g\longrightarrow \gl(\g)$ by
$$\ad_{x_1,x_2}y=[x_1,x_2,y].$$
Then, the map $\ad$ defines a representation of the $3$-Lie algebra $\g$ on itself, which we call {\bf adjoint representation} of $\g$.
 The coboundary operator associated to this representation is denoted by $\delta_{\g}$.

 \begin{lem}
If $\pi:\wedge^3\g\longrightarrow\g$ is a $3$-Lie bracket, then we have
\begin{equation}
\nrn{\pi,\alpha}=\delta_{\g}(\alpha),\quad \forall\alpha\in C^p(\g,\g),p\geq0.
\end{equation}
\end{lem}

\section{Generalized  representations of $3$-Lie algebras }

In this section,    we introduce a concept of generalized representation of a 3-Lie algebra using canonical structures.  First, we show that a representation of a 3-Lie algebra will give rise to a canonical structure.

Let $\rho:\wedge^2\g\longrightarrow\gl(V)$ be a linear map. Then, it induces a linear map $\bar{\rho}:\wedge^3(\g\oplus V)\longrightarrow \g\oplus V$ defined by
\begin{equation}
  \bar{\rho}(x+u,y+v,z+w)=\rho(x,y)(w)+\rho(y,z)(u)+\rho(z,x)(v),\quad \forall x,y,z\in\g, u,v,w\in V.
\end{equation}
Consider the graded Lie algebra given in Theorem \ref{thm:gradelie} associated to the vector space $\g\oplus V$.
\begin{pro}
 A linear map $\rho:\wedge^2\g\longrightarrow\gl(V)$ is a representation of the $3$-Lie algebra $\g$ on $V$ if and only if  $\pi+\bar{\rho}$ is a canonical structure in the graded Lie algebra associated to   $\g\oplus V$, i.e.
 $$
  \nrn{\pi+\bar{\rho},\pi+\bar{\rho}}=0.
 $$
\end{pro}

\pf By Lemma \ref{lem:semidirectp}, $\rho:\wedge^2\g\longrightarrow\gl(V)$ is a representation of $\g$ if and only if $\g\oplus V$ is a 3-Lie algebra, where the 3-Lie bracket is exactly given by
 \begin{eqnarray*}
[x+u,y+v,z+w]_{\rho}&=&[x,y,z]+\rho(x,y)(w)+\rho(y,z)(u)+\rho(z,x)(v)\\
&=&(\pi+\bar{\rho})(x+u,y+v,z+w).
\end{eqnarray*}
Thus, by Lemma \ref{lem:cs},  $\rho:\wedge^2\g\longrightarrow\gl(V)$ is a representation of $\g$ if and only if $\pi+\bar{\rho}$ is a canonical structure. \qed\vspace{3mm}

Now, we define  a new concept of representation of a 3-Lie algebra, generalizing the usual one. A so called generalized representation of a 3-Lie algebra involves two linear maps.
\begin{defi}\label{representation new}
 A {\bf generalized representation } of a $3$-Lie algebra $\frkg$ on a vector space $V$ consists of linear maps $\rho:\wedge^2 \frkg   \longrightarrow \gl(V)$ and $\nu:  \g\longrightarrow\Hom(\wedge^2 V, V)$ such that
 \begin{equation}
   \nrn{\pi+\bar{\rho}+\bar{\nu},\pi+\bar{\rho}+\bar{\nu}}=0,
 \end{equation}
 where $\bar{\nu}:\wedge^3( \g\oplus V)\longrightarrow \g\oplus V$ is induced by $\nu$ via
 $$
 \bar{\nu}(x+u,y+v,z+w)=\nu(x)(v\wedge w)+\nu(y)(w\wedge u)+\nu(z)(u\wedge v),\quad \forall~x,y,z\in\g, u,v,w\in V.
 $$
\end{defi}
We will refer to  a generalized representation by a triple  $(V;\rho,\nu).$

\begin{rmk}
  If $\nu=0$, then we recover the usual definition of a representation of a $3$-Lie algebra on a vector space (in the sense of Kasymov (\cite{Kasymov})). If the dimension of the vector space $V$ is $1$, then $\nu$ must be zero. In this case, we only have the usual representation.
\end{rmk}

 Given linear maps $\rho:\wedge^2 \frkg \longrightarrow \End(V)$ and $\nu:\frkg \longrightarrow \Hom(\wedge^2 V,V),$ define a trilinear bracket operation on $\g\oplus V$ by
\begin{eqnarray}
\nonumber[x+u,y+v,z+w]_{(\rho,\nu)}&=&[x,y,z]+\rho(x,y)(w)+\rho(y,z)(u)+\rho(z,x)(v)\\
\label{eq:formula}&&+\nu(x)(v\wedge w)+\nu(y)(w\wedge u)+\nu(z)(u\wedge v).
\end{eqnarray}

\begin{thm}\label{thm:semidirectproduct}
Let $(\g,[\cdot,\cdot,\cdot])$ be a $3$-Lie algebra and $(V;\rho,\nu)$ a generalized representation of $\g$. Then  $(\g\oplus V,[\cdot,\cdot,\cdot]_{(\rho,\nu)})$ is a $3$-Lie algebra, where $[\cdot,\cdot,\cdot]_{(\rho,\nu)}$ is given by \eqref{eq:formula}.
\end{thm}
We call the 3-Lie algebra $(\g\oplus V,[\cdot,\cdot,\cdot]_{(\rho,\nu)})$  the {\bf generalized semidirect product} of $\g$ and $V$.

\pf It follows  from
\begin{eqnarray*}
\nonumber[x+u,y+v,z+w]_{(\rho,\nu)}=(\pi+\bar{\rho}+\bar{\nu})(x+u,y+v,z+w)
\end{eqnarray*}
and Lemma \ref{lem:cs}.\qed\vspace{3mm}

In the following, we give a characterization of a generalized representation of a 3-Lie algebra.

\begin{pro}
Linear maps $\rho:\wedge^2 \frkg \longrightarrow \End(V)$ and $\nu:\frkg \longrightarrow \Hom(\wedge^2 V,V)$ give rise to a generalized representation of a   $3$-Lie algebra $\frkg$ on a vector space $V$ if and only if for all $x_i\in\frkg$, $v_j\in V$, the following equalities hold:
\begin{eqnarray}
\label{eq:r1}\rho(x_1,x_2)\rho(x_3,x_4)&=&\rho([x_1,x_2,x_3],x_4)+\rho(x_3,[x_1,x_2,x_4])+\rho(x_3,x_4)\rho(x_1,x_2),\\
\label{eq:r2}\rho(x_1,[x_2,x_3,x_4])&=&\rho(x_3,x_4)\rho(x_1,x_2)-\rho(x_2,x_4)\rho(x_1,x_3)+\rho(x_2,x_3)\rho(x_1,x_4),\\
\nonumber\rho(x_1,x_2)(\nu(x_3)(v_1,v_2))&=&\nu([x_1,x_2,x_3])(v_1,v_2)+\nu(x_3)(\rho(x_1,x_2)(v_1),v_2)\\
\label{eq:r3}&&+\nu(x_3)(v_1,\rho(x_1,x_2)(v_2)),\\
\nonumber\nu(x_1)(v_1,\rho(x_2,x_3)(v_2))&=&\nu(x_3)(v_2,\rho(x_2,x_1)(v_1))+\nu(x_2)(\rho(x_3,x_1)(v_1),v_2)\\
\label{eq:r4}&&+\rho(x_2,x_3)(\nu(x_1)(v_1,v_2)),\\
\label{eq:r6}\nu(x_1)(v_1,\nu(x_2)(v_2,v_3))&=&\nu(x_2)(\nu(x_1)(v_1,v_2),v_3)+\nu(x_2)(v_2,\nu(x_1)(v_1,v_3)),\\
\label{eq:r7}\nu(x_1)(\nu(x_2)(v_1,v_2),v_3)&=&\nu(x_2)(\nu(x_1)(v_1,v_2),v_3).
\end{eqnarray}
\end{pro}
\pf The triple $(V;\rho,\nu)$ is a generalized representation if and only if  $\nrn{\pi+\bar{\rho}+\bar{\nu},\pi+\bar{\rho}+\bar{\nu}}=0$. By straightforward computations,
$$
\nrn{\pi+\bar{\rho}+\bar{\nu},\pi+\bar{\rho}+\bar{\nu}}(x_1,x_2,x_3,x_4,v)=0
$$
is equivalent to \eqref{eq:r1}; And
$$
\nrn{\pi+\bar{\rho}+\bar{\nu},\pi+\bar{\rho}+\bar{\nu}}(x_1,v,x_2,x_3,x_4)=0
$$
is equivalent to \eqref{eq:r2}. Other identities  can be proved similarly.  The details are omited. \qed

\begin{rmk}
  By \eqref{eq:r1} and \eqref{eq:r2}, the map $\rho$ in a generalized representation $(V;\rho,\nu)$ gives rise to a usual representation in the sense of Definition \ref{defi:usualrep}. Conversely, for any representation $\rho$, $(V;\rho,\nu=0)$ is a generalized representation.
\end{rmk}
\begin{rmk}
Eq. \eqref{eq:r3} may be written as 
$\nu([x_1,x_2,x_3])=[\rho(x_1,x_2),\nu(x_3)]_{NR}, $  where $[\cdot,\cdot]_{NR}$ is the Nijenhuis-Richardson bracket on the graded vector space $\oplus_k\Hom(\wedge^kV,V)$; see \cite{NR} for more details about Nijenhuis-Richardson bracket.
\end{rmk}
\begin{rmk}
One may obtain from previous conditions the following identity
\begin{eqnarray}
\nonumber\nu([x_1,x_2,x_3])(v_1,v_2)&=&\rho(x_2,x_3)(\nu(x_1)(v_1,v_2))+\rho(x_3,x_1)(\nu(x_2)(v_1,v_2))\\
\label{eq:r5}&&+\rho(x_1,x_2)(\nu(x_3)(v_1,v_2)).
\end{eqnarray}
\end{rmk}
In the following, we provide a series of examples to illustrate the new concept of generalized representation.
\begin{ex}{\rm
 Let $\g$ be an abelian 3-Lie algebra. Define $\rho=0$, and $\nu=\xi\otimes \pi$, where  $\xi\in\g^*$ and $\pi\in\Hom(\wedge^2V\otimes V)$ is a Lie algebra structure on $V$. Then $(V;\rho,\nu)$ is a generalized representation. In fact,
 since $\g$ is abelian and $\rho=0$, \eqref{eq:r1}-\eqref{eq:r4} hold naturally.  Since $\pi$ satisfies the Jacobi identity, \eqref{eq:r6} and \eqref{eq:r7} also hold.
 }
\end{ex}

\begin{ex}{\rm
  Let $\g$ be the 3-dimensional 3-Lie algebra defined with respect to a  basis $\{x_1,x_2,x_3\}$ by the skew-symmetric bracket $[x_1,x_2,x_3]=x_1$.

 Let $V$ be a 2-dimensional vector space and $\{v_1,v_2\}$ its basis. The following  map $\rho$,  given with respect to bases by
 \begin{align*}
&  \rho(x_1,x_2)(v_1)=0,  \;
\rho(x_1,x_2)(v_2)=r_1 v_1 ,   \;
 \rho(x_1,x_3)(v_1)=0,\\
 & \rho(x_1,x_3)(v_2)=r_2 v_1, \;
 \rho(x_2,x_3)(v_1)=r_3 v_1,    \;  \rho(x_2,x_3)(v_2)=r_4 v_1+(r_3-1)v_2 ,
 \end{align*}
 where $r_1,\cdots,r_4$ are parameters in $\K$, 
defines a  representation of $\g$ on $V$ (in the usual sense, \eqref{eq:r1} and \eqref{eq:r2} are satisfied).

 It defines together with a map  $\nu$ a generalized representation of $\g$ on $V$ in the following cases.

\begin{enumerate}
\item
$   \rho(x_1,x_2)(v_1)=0,  \;
\rho(x_1,x_2)(v_2)=r_1 v_1 ,   \;
 \rho(x_1,x_3)(v_1)=0,$ \\
 $ \rho(x_1,x_3)(v_2)=r_2 v_1, \;
 \rho(x_2,x_3)(v_1)= v_1,    \;  \rho(x_2,x_3)(v_2)=r_4 v_1,$
 \\
$ \nu(x_1)(v_1,v_2)=0,  \;
\nu(x_2)(v_1,v_2)=s_1v_1,  \;
\nu(x_3)(v_1,v_2)=\frac{s_1 r_2}{r_1}v_1, $
\item
$   \rho(x_1,x_2)(v_1)=0,  \;
\rho(x_1,x_2)(v_2)=0 ,   \;
 \rho(x_1,x_3)(v_1)=0,$ \\
 $ \rho(x_1,x_3)(v_2)=r_2 v_1, \;
 \rho(x_2,x_3)(v_1)=v_1,    \;  \rho(x_2,x_3)(v_2)=r_4v_2 ,$
 \\
$ \nu(x_1)(v_1,v_2)=0,  \;
\nu(x_2)(v_1,v_2)=0,  \;
\nu(x_3)(v_1,v_2)=s_1 v_1, $
\item
$   \rho(x_1,x_2)(v_1)=0,  \;
\rho(x_1,x_2)(v_2)=0 ,   \;
 \rho(x_1,x_3)(v_1)=0,$ \\
 $ \rho(x_1,x_3)(v_2)=0, \;
 \rho(x_2,x_3)(v_1)=0,    \;  \rho(x_2,x_3)(v_2)=r_4 v_1-v_2 ,$
 \\
$ \nu(x_1)(v_1,v_2)=0,  \;
\nu(x_2)(v_1,v_2)=-s_1 r_4v_1+s_1 v_2,  \;
\nu(x_3)(v_1,v_2)=-s_2 r_4  v_1+s_2 v_2,$
\item
$   \rho(x_1,x_2)(v_1)=0,  \;
\rho(x_1,x_2)(v_2)=0 ,   \;
 \rho(x_1,x_3)(v_1)=0,$ \\
 $ \rho(x_1,x_3)(v_2)=0, \;
 \rho(x_2,x_3)(v_1)= v_1,    \;  \rho(x_2,x_3)(v_2)=r_4 v_1,$
 \\
$ \nu(x_1)(v_1,v_2)=0,  \;
\nu(x_2)(v_1,v_2)=s_1 v_1,  \;
\nu(x_3)(v_1,v_2)=s_2 v_1,$
\end{enumerate}
where $r_i$ and $s_j$ are parameters in $\K$.
}
\end{ex}

\begin{ex}{\rm
We may also consider, for the previous 3-dimensional 3-Lie algebra, the following map $\rho$  which defines  a  representation of $\g$ on $V$ (in the usual sense).
 \begin{align*}
&  \rho(x_1,x_2)(v_1)=r_1r_4v_1+r_1 v_2,  \;
\rho(x_1,x_2)(v_2)=-r_1r_4^2  v_1-r_1r_4 v_2 ,   \;
 \rho(x_1,x_3)(v_1)=r_2r_4v_1+r_2v_2,\\
 & \rho(x_1,x_3)(v_2)=-r_2 r_4^2 v_1-r_2 r_4 v_2, \;
 \rho(x_2,x_3)(v_1)=r_3 v_1,    \;  \rho(x_2,x_3)(v_2)=r_4 v_1+(1+r_3)v_2 .
 \end{align*}

 It defines together with a map  $\nu$ a generalized representation of $\g$ on $V$ in the following cases.

\begin{enumerate}
\item
$   \rho(x_1,x_2)(v_1)=r_1r_4v_1+r_1 v_2,  \;
\rho(x_1,x_2)(v_2)=-r_1r_4^2  v_1-r_1r_4 v_2 ,   \;
 \rho(x_1,x_3)(v_1)=r_2r_4v_1+r_2v_2,$ \\
 $  \rho(x_1,x_3)(v_2)=-r_2 r_4^2 v_1-r_2 r_4 v_2, \;
 \rho(x_2,x_3)(v_1)=0,    \;  \rho(x_2,x_3)(v_2)=r_4 v_1+v_2  ,$
 \\
$ \nu(x_1)(v_1,v_2)=0,  \;
\nu(x_2)(v_1,v_2)=s_1 r_4v_1+s_1 v_2,  \;
\nu(x_3)(v_1,v_2)=\frac{s_1 r_2 r_4}{r_1}v_1+\frac{s_1 r_2}{r_1} v_2, $
\item
$   \rho(x_1,x_2)(v_1)=0,  \;
\rho(x_1,x_2)(v_2)=0 ,   \;
 \rho(x_1,x_3)(v_1)=r_2r_4v_1+r_2v_2,$ \\
 $  \rho(x_1,x_3)(v_2)=-r_2 r_4^2 v_1-r_2 r_4 v_2, \;
 \rho(x_2,x_3)(v_1)=0,    \;  \rho(x_2,x_3)(v_2)=r_4 v_1+v_2  ,$
 \\
$ \nu(x_1)(v_1,v_2)=0,  \;
\nu(x_2)(v_1,v_2)=0,  \;
\nu(x_3)(v_1,v_2)=r_4 s_1v_1+s_1 v_2, $
\item
$   \rho(x_1,x_2)(v_1)=0  \;
\rho(x_1,x_2)(v_2)=0 ,   \;
 \rho(x_1,x_3)(v_1)=0,$ \\
 $  \rho(x_1,x_3)(v_2)=0, \;
 \rho(x_2,x_3)(v_1)=0,    \;  \rho(x_2,x_3)(v_2)=r_4 v_1+v_2 ,$
 \\
$ \nu(x_1)(v_1,v_2)=0,  \;
\nu(x_2)(v_1,v_2)=s_1 r_4v_1+s_1 v_2,  \;
\nu(x_3)(v_1,v_2)=s_2 r_4v_1+s_2 v_2, $
\item
$   \rho(x_1,x_2)(v_1)=0,  \;
\rho(x_1,x_2)(v_2)=0 ,   \;
 \rho(x_1,x_3)(v_1)=0,$ \\
 $  \rho(x_1,x_3)(v_2)=0, \;
 \rho(x_2,x_3)(v_1)=- v_1,    \;  \rho(x_2,x_3)(v_2)=r_4 v_1 ,$
 \\
$ \nu(x_1)(v_1,v_2)=0,  \;
\nu(x_2)(v_1,v_2)=s_1v_1,  \;
\nu(x_3)(v_1,v_2)=s_2v_1, $
\item
$   \rho(x_1,x_2)(v_1)=0,  \;
\rho(x_1,x_2)(v_2)=0 ,   \;
 \rho(x_1,x_3)(v_1)=r_2v_2,$ \\
 $  \rho(x_1,x_3)(v_2)=0, \;
 \rho(x_2,x_3)(v_1)=0,    \;  \rho(x_2,x_3)(v_2)=v_2  ,$
 \\
$ \nu(x_1)(v_1,v_2)=0,  \;
\nu(x_2)(v_1,v_2)=0,  \;
\nu(x_3)(v_1,v_2)=s_1 v_2, $

\item
$   \rho(x_1,x_2)(v_1)=r_1 v_2,  \;
\rho(x_1,x_2)(v_2)=0 ,   \;
 \rho(x_1,x_3)(v_1)=0,$ \\
 $  \rho(x_1,x_3)(v_2)=0, \;
 \rho(x_2,x_3)(v_1)=0,    \;  \rho(x_2,x_3)(v_2)=v_2  ,$
 \\
$ \nu(x_1)(v_1,v_2)=0,  \;
\nu(x_2)(v_1,v_2)=s_1 v_2,  \;
\nu(x_3)(v_1,v_2)=0,$

\end{enumerate}
where $r_i$ and $s_j$ are parameters in $\K$.
}
\end{ex}

In the following, we consider generalized representations where the usual representation vanishes.
\begin{pro}
Let $\g$ be  an $n$-dimensional  $3$-Lie algebra such that $[\g,\g,\g]=\g$ ($\g$ is said to be perfect). If the generalized  representation $(V;\rho,\nu)$ is such that $\rho =0$ then $\nu =0$.
\end{pro}
\pf
It follows from condition \eqref{eq:r3}.\qed

\begin{cor}
Let $\g$ be  the $4$-dimensional simple $3$-Lie algebra defined with respect to a  basis $\{x_1,x_2,x_3,x_4\}$ by the skew-symmetric brackets
$$[x_1,x_2,x_3]=x_4, \;[x_1,x_2,x_4]=x_3, \; [x_1,x_3,x_4]=x_2,\; [x_2,x_3,x_4]=x_1. $$
If an $n$-dimensional  generalized representation $(V;\rho,\nu)$ of  $\g$ satisfies $\rho=0$, then $\nu=0$.
\end{cor}
Now, we show some examples where $[\g,\g,\g]\neq \g$.

\begin{ex}{\rm
 We provide for the 3-dimensional  3-Lie algebra, defined with respect to a  basis $\{x_1,x_2,x_3\}$ by the skew-symmetric bracket $[x_1,x_2,x_3]=x_1$, all the 2-dimensional generalized representations $(V;\rho,\nu)$  with $\rho=0$.
 Every generalized representation  on  a 2-dimensional vector space $V$ with a trivial $\rho$ is given by one of the following maps $\nu$ defined, with respect to a basis  $\{ v_1,v_2\}$ of $V$, by

\begin{enumerate}
\item
$ \nu(x_1)(v_1,v_2)=0,  \;
\nu(x_2)(v_1,v_2)=s_1v_1+s_2 v_2,  \;
\nu(x_3)(v_1,v_2)=s_3 v_1+\frac{s_2 s_3}{s_1}v_2, $
\item
$ \nu(x_1)(v_1,v_2)=0,  \;
\nu(x_2)(v_1,v_2)=s_1v_2,  \;
\nu(x_3)(v_1,v_2)=s_2 v_2, $
\item
$ \nu(x_1)(v_1,v_2)=0,  \;
\nu(x_2)(v_1,v_2)=0,  \;
\nu(x_3)(v_1,v_2)=s_1 v_1+s_2 v_2,$
\end{enumerate}
where $s_1,s_2$ are parameters in $\K$.
}
\end{ex}

\begin{ex}{\rm
Let $\g$ be  the $4$-dimensional  $3$-Lie algebra defined, with respect to a  basis $\{x_1,x_2,x_3,\\x_4\}$, by the skew-symmetric brackets
$$[x_1,x_2,x_4]=x_3,\;[x_1,x_3,x_4]=x_2, \; [x_2,x_3,x_4]=x_1. $$
Every generalized representation $(V;\rho,\nu)$, on a 2-dimensional vector space $V$  with a trivial $\rho$, of $\g$ is given by one of the following maps $\nu$, defined with respect to a basis  $\{ v_1,v_2\}$ of $V$, by
$$ \nu(x_1)(v_1,v_2)=0,  \;
\nu(x_2)(v_1,v_2)=0,  \;
\nu(x_3)(v_1,v_2)=0,  \;
\nu(x_4)(v_1,v_2)=s_1 v_1+s_2 v_2, $$
where $s_1,s_2$ are parameters in $\K$.
}
\end{ex}

\begin{ex}{\rm
Let $\g$ be  the 4-dimensional  3-Lie algebra defined, with respect to a basis $\{x_1,x_2,x_3,x_4\}$,  by $ [x_2,x_3,x_4]=x_1$.
 Every generalized representation $(V;\rho,\nu)$, on a 2-dimensional vector space  $V$  with trivial $\rho$, of $\g$ is given by one of the following maps $\nu$ defined, with respect to a basis  $\{ v_1,v_2\}$ of $V$, by

\begin{enumerate}\item $  \nu(x_1)(v_1,v_2)=0,  \;
\nu(x_2)(v_1,v_2)=s_1 v_1+s_2 v_2,  \; \\ $ $
\nu(x_3)(v_1,v_2)=s_3 v_1+\frac{s_2 s_3}{s_1} v_2, \;
\nu(x_4)(v_1,v_2)=s_4 v_1+\frac{s_2 s_4}{s_1} v_2, $
\item  $ \nu(x_1)(v_1,v_2)=0,  \;
\nu(x_2)(v_1,v_2)=0,  \; \\ $ $
\nu(x_3)(v_1,v_2)=s_1 v_1+s_2 v_2, \;
\nu(x_4)(v_1,v_2)=s_3 v_1+\frac{s_2 s_4}{s_1} v_2, $
\end{enumerate}
where $s_1,s_2,s_3,s_4$ are parameters in $\K$.
}
\end{ex}

Now, we characterize the notion of equivalent generalized representations of 3-Lie algebras.

\begin{defi}
  Let $(V_1;\rho_1,\nu_1)$ and $(V_2;\rho_2,\nu_2)$ be two generalized representations of a $3$-Lie algebra $(\g,[\cdot,\cdot,\cdot])$. They are said to be {\bf equivalent} if there exists an isomorphism of vector spaces $T:V_1\longrightarrow V_2$ such that
  $$
  T\rho_1(x,y)(u)=\rho_2(x,y)(Tu),\quad T\nu_1(x)(u,v)=\nu_2(x)(Tu,Tv),\quad\forall x,y\in\g,~u, v\in V_1.
  $$
  In terms of diagrams, we have
  $$
\xymatrix{
 \wedge^2\g\times V_1 \ar[d]_{ id\times T }\ar[rr]^{\rho_1}
                && V_1  \ar[d]^{T}  \\
 \wedge^2\g\times V_2 \ar[rr]^{\rho_2}
                && V_2  },\quad \xymatrix{
 \g\times \wedge^2 V_1 \ar[d]_{ id\times \wedge^2T }\ar[rr]^{\nu_1}
                && V_1  \ar[d]^{T}  \\
 \g\times \wedge^2 V_2 \ar[rr]^{\nu_2}
                && V_2.  }
$$
\end{defi}

\section{New cohomology complex of $3$-Lie algebras }
Based on the generalized representations defined in the previous section, we introduce a new type of cohomology for 3-Lie algebras.

Let $(\g, [\cdot,\cdot,\cdot ])$ be a 3-Lie algebra and $(V;\rho,\nu)$ a generalized representation of $\g$.
We set $C^p_>(\g\oplus V,V)$ to be the set of $(p+1)$-cochains, which is defined as a subset of $C^p(\g\oplus V,V)$ such that
\begin{equation}\label{eq:C}
C^p(\g\oplus V,V)={C}^p_>(\g\oplus V,V)\oplus C^p(V,V).
\end{equation}
Recall that $C^p(\g\oplus V,V)=\{ \alpha:\wedge^2(\frkg\oplus V)\otimes\stackrel{(p \text{ times})}{\cdots }\otimes\wedge^2(\frkg\oplus V)\wedge(\frkg\oplus V)\longrightarrow V\}.$

By direct calculation, we have
$$\nrn{\pi+\bar{\rho}+\bar{\nu}, {C}^\bullet_>(\g\oplus V,V)}\subseteq{C}^{\bullet+1}_>(\g\oplus V,V).$$
 Define $\dM:{C}^p_>(\g\oplus V,V)\longrightarrow{C}^{p+1}_>(\g\oplus V,V)$  by
\begin{equation}
\dM(\alpha):=\nrn{\pi+\bar{\rho}+\bar{\nu},\alpha},\quad \alpha\in {C}^p_>(\g\oplus V,V).
\end{equation}

\begin{thm}\label{thm:dd0}
Let $(V;\rho,\nu)$ be a generalized representation of a $3$-Lie algebra $\g$. Then $\dM\circ \dM =0.$ Thus, we obtain a new cohomology complex, where the space of $p$-cochains is given by ${C}^{p-1}_>(\g\oplus V,V)$.
\end{thm}
\pf By the graded Jacobi identity, for any $\alpha\in{C}^{p-1}_>(\g\oplus V,V)$, we have
$$
\dM\circ \dM(\alpha):=\nrn{\pi+\bar{\rho}+\bar{\nu},\nrn{\pi+\bar{\rho}+\bar{\nu},\alpha}}=\half \nrn{\nrn{\pi+\bar{\rho}+\bar{\nu},\pi+\bar{\rho}+\bar{\nu}},\alpha}=0,
$$
which finishes the proof.\qed\vspace{3mm}

An element $\alpha\in {C}^{p-1}_>(\g\oplus V,V)$ is called a $p$-cocycle if $\dM(\alpha)=0$; It is called a $p$-coboundary if there exists  $\beta\in {C}^{p-2}_>(\g\oplus V,V)$ such that $\alpha=\dM(\beta)$. Denote by $\mathcal{Z}^p(\g;V)$ and $\mathcal{B}^p(\g;V)$ the sets of $p$-cocycles and  $p$-coboundaries respectively. By Theorem \ref{thm:dd0}, we have $\mathcal{B}^p(\g;V)\subset\mathcal{Z}^p(\g;V)$. We define the $p$-th cohomolgy group 
$\mathcal{H}^p(\g;V)$ to be $\mathcal{Z}^p(\g;V)/\mathcal{B}^p(\g;V)$.\\

A relationship  between this new cohomology and the one given by \eqref{eq:cohomology}, is stated in the following result. 

\begin{pro}
  There is  a forgetful map from $\mathcal{H}^p(\g;V)$ to $H^p(\g;V)$.
\end{pro}
\pf
It is obvious that ${C}^p(\g,V)\subseteq{C}^p_>(\g\oplus V,V)$. By direct calculation, for $\frkX_i\in\wedge^2\g,z\in\g$, we have
$$\dM(\alpha)(\frkX_1,\cdots ,\frkX_{p+1},z)=\delta_\rho(\alpha)(\frkX_1,\cdots ,\frkX_{p+1},z),\quad \alpha\in {C}^p(\g,V),$$
where $\delta_\rho$ is the coboundary operator given by \eqref{eq:drho}. Thus, the natural projection from ${C}^p_>(\g\oplus V,V)$ to ${C}^p(\g,V)$ induces a forgetful map from $\mathcal{H}^p(\g;V)$ to $H^p(\g;V)$. \qed\vspace{3mm}

In the sequel, we give some characterization of  low dimensional cocycles.

\begin{pro}
A linear map $\alpha\in\Hom(\g,V)$ is a $1$-cocycle if only if for all $x_1,x_2,x_3\in \g,v\in V$, the following identities hold :
\begin{eqnarray*}
  & \nu(x_1)(\alpha(x_2),v)-\nu(x_2)(\alpha(x_1),v)=0,\\
  & \alpha([x_1,x_2,x_3])-\rho(x_1,x_2)(\alpha(x_3))-\rho(x_2,x_3)(\alpha(x_1))-\rho(x_3,x_1)(\alpha(x_2))=0.
\end{eqnarray*}
\end{pro}
\pf
For $\alpha\in\Hom(\g,V)$, we have
\begin{eqnarray*}
\dM(\alpha)(x_1,x_2,v)=\nu(x_1)(\alpha(x_2),v)-\nu(x_2)(\alpha(x_1),v),
\end{eqnarray*}
and
\begin{eqnarray*}
\dM(\alpha)(x_1,x_2,x_3)&=&\delta_\rho(\alpha)(x_1,x_2,x_3)\\
&=&\rho(x_1,x_2)(\alpha(x_3))+\rho(x_2,x_3)(\alpha(x_1))+\rho(x_3,x_1)(\alpha(x_2))-\alpha([x_1,x_2,x_3]), \end{eqnarray*}
which finishes the proof. \qed

\begin{pro}
  A $2$-cochain $\alpha_1+\alpha_2+\alpha_3\in {C}^{1}_>(\g\oplus V,V)$, where $ \alpha_1\in\Hom(\wedge^2 V\wedge \g,V),~\alpha_2\in\Hom(\wedge^2\g\wedge V,V),~\alpha_3\in\Hom(\wedge^3\g,V)$, is a $2$-cocycle if and only if for all $x_i\in\g,v_j\in V$ and $v\in V$, the following identities hold:
  \begin{eqnarray}
  \nonumber 0&=&-\rho(x_1,x_2)(\alpha_3(x_3,x_4,x_5))-\alpha_3(x_1,x_2,[x_3,x_4,x_5])+\rho(x_4,x_5)(\alpha_3(x_1,x_2,x_3))
 \\&& \nonumber+\alpha_3([x_1,x_2,x_3],x_4,x_5)+\rho(x_5,x_3)(\alpha_3(x_1,x_2,x_4))+\alpha_3(x_3,[x_1,x_2,x_4],x_5)
 \\\label{eq:2cocycle1}&&+\rho(x_3,x_4)(\alpha_3(x_1,x_2,x_5))+\alpha_3(x_3,x_4,[x_1,x_2,x_5]),\\
\nonumber0&=&\nu(x_4)(\alpha_3(x_1,x_2,x_3),v)+\nu(x_3)(v,\alpha_3(x_1,x_2,x_4))+\rho(x_1,x_2)(\alpha_2(x_3,x_4,v))
\\&&\label{eq:2cocycle2}-\rho(x_3,x_4)(\alpha_2(x_1,x_2,v))-\alpha_2([x_1,x_2,x_3],x_4,v)-\alpha_2(x_3,[x_1,x_2,x_4],v),
\\
\nonumber0&=&\nu(x_1)(v,\alpha_3(x_2,x_3,x_4))+\rho(x_3,x_4)(\alpha_2(x_1,x_2,v))-\rho(x_2,x_4)(\alpha_2(x_1,x_3,v))
\\&&\nonumber+\rho(x_2,x_3)(\alpha_2(x_1,x_4,v))+\alpha_2(x_3,x_4,\rho(x_1,x_2)(v))-\alpha_2(x_2,x_4,\rho(x_1,x_3)(v))
\\&&\label{eq:2cocycle3}+\alpha_2(x_2,x_3,\rho(x_1,x_4)(v))-\alpha_2(x_1,[x_2,x_3,x_4],v),\\
\nonumber0&=&\nu(x_3)(v_2,\alpha_2(x_1,x_2,v_1))+\nu(x_3)(\alpha_2(x_1,x_2,v_2),v_1)+\alpha_2(x_1,x_2,\nu(x_3)(v_1,v_2))
\\&&\nonumber+\rho(x_1,x_2)(\alpha_1(v_1,v_2,x_3))-\alpha_1(\rho(x_1,x_2)(v_1),v_2,x_3)-\alpha_1(v_1,\rho(x_1,x_2)(v_2),x_3)
\\&&\label{eq:2cocycle4}-\alpha_1(v_1,v_2,[x_1,x_2,x_3]),\\
\nonumber0&=&\nu(x_3)(v_2,\alpha_2(x_2,x_1,v_1))+\nu(x_2)(\alpha_2(x_3,x_1,v_1),v_2)-\nu(x_1)(v_1,\alpha_2(x_2,x_3,v_2))
\nonumber\\&&\nonumber+\alpha_2(x_2,x_3,\nu(x_1)(v_1,v_2))+\rho(x_2,x_3)(\alpha_1(v_1,v_2,x_1))+\alpha_1(\rho(x_1,x_2)(v_1),v_2,x_3)
\\&&\label{eq:2cocycle5}-\alpha_1(v_1,\rho(x_2,x_3)(v_2),x_1)+\alpha_1(v_2,\rho(x_1,x_3)(v_1),x_2),\\
\nonumber0&=&\alpha_2(x_1,x_3,\nu(x_2)(v_1,v_2))-\alpha_2(x_2,x_3,\nu(x_1)(v_1,v_2))-\alpha_2(x_1,x_2,\nu(x_3)(v_1,v_2))
\\&&\nonumber+\rho(x_1,x_3)(\alpha_1(v_1,v_2,x_2))-\rho(x_1,x_2)(\alpha_1(v_1,v_2,x_3))-\rho(x_2,x_3)(\alpha_1(v_1,v_2,x_1))
\\&&\label{eq:2cocycle6}+\alpha_1(v_1,v_2,[x_1,x_2,x_3]),\\
\nonumber0&=&-\nu(x_2)(\alpha_1(v_1,v_2,x_1),v_3)-\nu(x_2)( v_2,\alpha_1(v_1,v_3,x_1))+\nu(x_1)( v_1,\alpha_1(v_2,v_3,x_2))
\\&&\label{eq:2cocycle7}-\alpha_1(\nu(x_1)(v_1,v_2),v_3,x_2)-\alpha_1(v_2,\nu(x_1)(v_1,v_3),x_2)+\alpha_1(v_1,\nu(x_2)(v_2,v_3),x_1),\\
\nonumber0&=&\nu(x_2)(\alpha_1(v_1,v_2,x_1), v_3)-\nu(x_1)( v_3,\alpha_1(v_1,v_2,x_2))+\alpha_1(\nu(x_1)(v_1,v_2),v_3,x_2)\\
&&\label{eq:2cocycle8}-\alpha_1(\nu(x_2)(v_1,v_2),v_3,x_1).
  \end{eqnarray}
\end{pro}
\pf
For $\alpha_3\in\Hom(\wedge^3\g,V)$, we have
\begin{eqnarray*}
\dM(\alpha_3)(x_1,x_2,x_3,x_4,x_5)&=&\rho(x_1,x_2)(\omega(x_3,x_4,x_5))+\omega(x_1,x_2,[x_3,x_4,x_5])\\&& \nonumber-\rho(x_4,x_5)(\omega(x_1,x_2,x_3))
 -\omega([x_1,x_2,x_3],x_4,x_5)\\&& \nonumber-\rho(x_5,x_3)(\omega(x_1,x_2,x_4))-\omega(x_3,[x_1,x_2,x_4],x_5)
 \\&& \nonumber-\rho(x_3,x_4)(\omega(x_1,x_2,x_5))-\omega(x_3,x_4,[x_1,x_2,x_5]),\\
\dM(\alpha_3)(x_1,x_2,x_3,x_4,v)&=&\nu(x_4)(\alpha_3(x_1,x_2,x_3),v)+\nu(x_3)(v,\alpha_3(x_1,x_2,x_4)),\\
\dM(\alpha_3)(x_1,v,x_2,x_3,x_4)&=&\nu(x_1)(v,\alpha_3(x_2,x_3,x_4)).
\end{eqnarray*}
For $\alpha_2\in\Hom(\wedge^2\g\wedge V,V)$, we have
\begin{eqnarray*}
\dM(\alpha_2)(x_1,x_2,x_3,x_4,v)&=&\rho(x_1,x_2)(\alpha_2(x_3,x_4,v))-\rho(x_3,x_4)(\alpha_2(x_1,x_2,v))\\
&&-\alpha_2([x_1,x_2,x_3],x_4,v)-\alpha_2(x_3,[x_1,x_2,x_4],v),\\
\dM(\alpha_2)(x_1,v,x_2,x_3,x_4)&=&\rho(x_3,x_4)(\alpha_2(x_1,x_2,v))-\rho(x_2,x_4)(\alpha_2(x_1,x_3,v))\\
&&+\rho(x_2,x_3)(\alpha_2(x_1,x_4,v))+\alpha_2(x_3,x_4,\rho(x_1,x_2)(v))\\
&&-\alpha_2(x_2,x_4,\rho(x_1,x_3)(v))+\alpha_2(x_2,x_3,\rho(x_1,x_4)(v))\\
&&-\alpha_2(x_1,[x_2,x_3,x_4],v),
\end{eqnarray*}
\begin{eqnarray*}
\dM(\alpha_2)(x_1,x_2,v_1,v_2,x_3)&=&\nu(x_3)(v_2,\alpha_2(x_1,x_2,v_1))+\nu(x_3)(\alpha_2(x_1,x_2,v_2),v_1)\\
&&+\alpha_2(x_1,x_2,\nu(x_3)(v_1,v_2)),\\
\dM(\alpha_2)(x_1,v_1,x_2,v_2,x_3)&=&\nu(x_3)(v_2,\alpha_2(x_2,x_1,v_1))+\nu(x_2)(\alpha_2(x_3,x_1,v_1),v_2)\\
&&-\nu(x_1)(v_1,\alpha_2(x_2,x_3,v_2))+\alpha_2(x_2,x_3,\nu(x_1)(v_1,v_2)),\\
\dM(\alpha_2)(v_1,v_2 ,x_1,x_2,x_3)&=&\alpha_2(x_1,x_3,\nu(x_2)(v_1,v_2))-\alpha_2(x_2,x_3,\nu(x_1)(v_1,v_2))\\&&-\alpha_2(x_1,x_2,\nu(x_3)(v_1,v_2)).
\end{eqnarray*}
For $\alpha_1\in\Hom(\wedge^2 V\wedge \g,V)$, we have
\begin{eqnarray*}
\dM(\alpha_1)(x_1,x_2,v_1,v_2,x_3)&=&\rho(x_1,x_2)(\alpha_1(v_1,v_2,x_3))-\alpha_1(\rho(x_1,x_2)(v_1),v_2,x_3)\\
&&-\alpha_1(v_1,\rho(x_1,x_2)(v_2),x_3)-\alpha_1(v_1,v_2,[x_1,x_2,x_3]),\\
\dM(\alpha_1)(x_1,v_1,x_2,v_2,x_3)&=&\rho(x_2,x_3)(\alpha_1(v_1,v_2,x_1))+\alpha_1(\rho(x_1,x_2)(v_1),v_2,x_3)\\
&&-\alpha_1(v_1,\rho(x_2,x_3)(v_2),x_1)+\alpha_1(v_2,\rho(x_1,x_3)(v_1),x_2),\\
\dM(\alpha_1)(v_1,v_2 ,x_1,x_2,x_3)&=&\rho(x_1,x_3)(\alpha_1(v_1,v_2,x_2))-\rho(x_1,x_2)(\alpha_1(v_1,v_2,x_3))\\
&&-\rho(x_2,x_3)(\alpha_1(v_1,v_2,x_1))+\alpha_1(v_1,v_2,[x_1,x_2,x_3]),\\
\dM(\alpha_1)(x_1,v_1,v_2,v_3,x_2)&=&-\nu(x_2)(\alpha_1(v_1,v_2,x_1),v_3)-\nu(x_2)( v_2,\alpha_1(v_1,v_3,x_1))\\&&+\nu(x_1)( v_1,\alpha_1(v_2,v_3,x_2))
-\alpha_1(\nu(x_1)(v_1,v_2),v_3,x_2)\\&&-\alpha_1(v_2,\nu(x_1)(v_1,v_3),x_2)+\alpha_1(v_1,\nu(x_2)(v_2,v_3),x_1),\\
\dM(\alpha_1)(v_1,v_2 ,x_1,x_2,v_3)&=&\nu(x_2)(\alpha_1(v_1,v_2,x_1), v_3)-\nu(x_1)( v_3,\alpha_1(v_1,v_2,x_2))\\
&&+\alpha_1(\nu(x_1)(v_1,v_2),v_3,x_2)-\alpha_1(\nu(x_2)(v_1,v_2),v_3,x_1).
\end{eqnarray*}
Thus, $\dM(\alpha_1+\alpha_2+\alpha_3)=0$ if and only if Eqs \eqref{eq:2cocycle1}-\eqref{eq:2cocycle8} hold. \qed\vspace{3mm}

In the following we provide two examples of computation of 2-cocycles of the 3-dimensional ternary algebra defined with respect to a basis $\{x_1,x_2,x_3\}$ by the bracket $[x_1,x_2,x_3]=x_1$. We consider two different generalized representations on a $2$-dimensional vector space $V$, with basis $\{v_1,v_2\}$.

\begin{ex}[Example 1]{\rm
We consider the generalized representation $(V,\rho,\nu)$, where   $\rho$ and $\nu$ are defined with respect to the basis by
\begin{eqnarray*}
& \rho (x_1,x_2)(v_1)=0,\
 \rho (x_1,x_2)(v_2)=0,\
 \rho (x_1,x_3)(v_1)=s_1 s_2 v_1+s_1 v_2,\\
 & \rho (x_1,x_3)(v_2)=-s_1 s_2^2 v_1-s_1 s_2 v_2,\
\rho (x_2,x_3)(v_1)=0,\
\rho (x_2,x_3)(v_2)=s_2 v_1+v_2,
\end{eqnarray*}
\begin{eqnarray*}
 \nu (x_1)(v_1,v_2)=0,\
 \nu (x_2)(v_1,v_2)=0,\
 \nu (x_3)(v_1,v_2)=s_3 s_2  v_1+s_3 v_2,
\end{eqnarray*}
with $s_1,s_2,s_3$   parameters in $\K$.

The 2-cocycles are given by $\alpha_1=0$, $\alpha_3=0$ and $\alpha_2$  defined as
 \begin{eqnarray*}
& \alpha_2 (x_1,x_2,v_1)=0,\
 & \alpha_2 (x_1,x_2,v_2)=0,\\
& \alpha_2 (x_1,x_3,v_1)=s_2 p_1 v_1+p_1v_2,\
 & \alpha_2 (x_1,x_3,v_2)=-s_2^2 p_1 v_1-s_2 p_1 v_2,\\
& \alpha_2 (x_2,x_3,v_1)=s_2 p_2 v_1+p_2 v_2,\
 & \alpha_2 (x_2,x_3,v_2)=-s_2^2 p_2 v_1-s_2 p_2 v_2,
\end{eqnarray*}
where $p_1,p_2$ are parameters in $\K$.
}
\end{ex}

\begin{ex}[Example 2]
{\rm
Now, we consider the generalized representation $(V,\rho,\nu)$, where $\rho$ is trivial and $\nu$ is given with respect to the basis by
\begin{eqnarray*}
 \nu (x_1)(v_1,v_2)=0,\
 \nu (x_2)(v_1,v_2)=s_1 v_1+s_2 v_2,\
 \nu (x_3)(v_1,v_2)=s_3 v_1+\frac{s_2 s_3}{s_1} v_2,
\end{eqnarray*}
with $s_1,s_2,s_3$  parameters in $\K$.

The 2-cocycles are given by $\alpha_1=0$ and $\alpha_2$, $\alpha_3$  defined as
 \begin{eqnarray*}
& \alpha_2 (x_1,x_2,v_1)=-s_1 p_2  v_1-s_2 p_2  v_2,\
 & \alpha_2 (x_1,x_2,v_2)=s_1 p_1 v_1+s_2 p_1 v_2,\\
& \alpha_2 (x_1,x_3,v_1)=-s_3 p_2  v_1- \frac{s_2 s_3 p_2}{s_1} v_2,\
 & \alpha_2 (x_1,x_3,v_2)=s_3 p_1 v_1+\frac{s_2 s_3  p_1}{s_1} v_2,\\
& \alpha_2 (x_2,x_3,v_1)=p_3 v_1+\frac{s_2 p_3}{s_1} v_2,\
 & \alpha_2 (x_2,x_3,v_2)=p_4 v_1+\frac{s_2 p_4}{s_1} v_2,
\end{eqnarray*}

\begin{equation*}
\alpha_3(x_1,x_2,x_3)=p_1 v_1+p_2 v_2,
\end{equation*}
where $p_1,p_2,p_3,p_4$ are parameters in $\K$.
}
\end{ex}

\section{Abelian extensions of $3$-Lie algebras }\label{sec:pre}

In this section, first we study the split abelian extension of 3-Lie algebras, which is isomorphic to a generalized semidirect product 3-Lie algebra. This provides a good motivation for the introduction of the  generalized representation. Then,  we study non-split abelian extensions of 3-Lie algebras. Unlike the case of Lie algebras, they cannot be described by 2-cocycles. Finally, we describe non-split abelian extensions via Maurer-Cartan elements.

Consider an exact sequence of   $3$-Lie algebras  :
\begin{equation}\label{eq:sequence}
\xymatrix@C=0.5cm{
  0 \ar[r] &V \ar[rr]^{i} && \hat{\frkg} \ar[rr]^{p} && \frkg \ar[r] & 0 },
  \end{equation}
then $\hat{\frkg}$ is said to be an abelian extension of the $3$-Lie algebra $\frkg$ by $V$ if $V$ is  abelian. A linear map  $\sigma:\frkg\longrightarrow\hat{\frkg}$ is called a splitting of $\hat{\g}$ if it satisfies $p\circ \sigma=id_{\frkg}$. If there exists a splitting which is also a homomorphism between 3-Lie algebras, we say that the abelian extension is {\bf split}.

Let $\hat{\g}$ be a split abelian extension and  $\sigma:\frkg\longrightarrow\hat{\frkg}$ the corresponding splitting.   Define $\rho:\wedge^2\g\longrightarrow\gl(V)$ and $\nu:\g\longrightarrow\Hom(\wedge^2V,V)$   by
\begin{eqnarray*}
 \rho(x,y)(u)&=&[\sigma(x),\sigma(y),u]_{\hat{\frkg}},\\
  \nu(x)(u,v)&=&[\sigma(x),u,v]_{\hat{\frkg}}.
\end{eqnarray*}
Then, we can transfer the $3$-Lie algebra  structure on $\hat{\frkg}$ to that on $\g\oplus V$ in terms of $\rho$ and $\nu$:
\begin{eqnarray*}
[x+u,y+v,z+w]_{(\rho,\nu)}&=&[x,y,z]+\rho(x,y)(w)+\rho(y,z)(u)+\rho(z,x)(v)\\
&&+\nu(x)(v\wedge w)+\nu(y)(w\wedge u)+\nu(z)(u\wedge v).
\end{eqnarray*}
Note that the Fundamental Identity  gives the character of  $\rho$ and $\nu$. However, by Theorem \ref{thm:semidirectproduct}, it is straightforward to obtain the following proposition.

\begin{pro}
  Any split abelian extension of $3$-Lie algebras is isomorphic to a generalized semidirect product $3$-Lie algebra.
\end{pro}

 Now, for non-split abelian extensions, we can further define $\omega:\wedge^3\frkg\longrightarrow V$ by
\begin{eqnarray*}
  \omega(x,y,z)&=&[\sigma(x),\sigma(y),\sigma(z)]_{\hat{\frkg}}-\sigma[x,y,z]_{\frkg}.
\end{eqnarray*}
Then, we also transfer the $3$-Lie algebra  structure on $\hat{\frkg}$ to that on $\g\oplus V$ in terms of $\rho,\nu$ and $\omega:$
\begin{eqnarray*}
[x_1+v_1,x_2+v_2,x_3+v_3]_{(\rho,\nu,\omega)}&=&[x_1,x_2,x_3]_{\frkg}+\rho(x_1,x_2)(v_3)+\rho(x_3,x_1)(v_2)+\rho(x_2,x_3)(v_1)\\
&&+\nu(x_1)(v_2,v_3)+\nu(x_2)(v_3,v_1)+\nu(x_3)(v_1,v_2)+\omega(x_1,x_2,x_3).
\end{eqnarray*}

The Fundamental Identity  gives the character of  $\rho$, $\nu$ and $\omega$.
\begin{thm}\label{thm:abelian ext}
 With above notations, $(\g\oplus V,[\cdot,\cdot,\cdot]_{(\rho,\nu,\omega)})$ is a $3$-Lie algebra if and only if for all $x_1,x_2,x_3,x_4,x_5\in\g$ and $v,v_1,v_2,v_3\in V$, Eqs. \eqref{eq:r3}-\eqref{eq:r7} and the following identities hold:
    \begin{eqnarray}
 \nonumber 0&=&-\rho(x_1,x_2)(\omega(x_3,x_4,x_5))-\omega(x_1,x_2,[x_3,x_4,x_5])+\rho(x_4,x_5)(\omega(x_1,x_2,x_3))
 \\&& \nonumber+\omega([x_1,x_2,x_3],x_4,x_5)+\rho(x_5,x_3)(\omega(x_1,x_2,x_4))+\omega(x_3,[x_1,x_2,x_4],x_5)
 \\\label{eq:t1}&&+\rho(x_3,x_4)(\omega(x_1,x_2,x_5))+\omega(x_3,x_4,[x_1,x_2,x_5]),\\
 \nonumber 0&=&\nu(x_1)(v,\omega(x_2,x_3,x_4))+\rho([x_2,x_3,x_4],x_1)(v)+\rho(x_3,x_4)\rho(x_1,x_2)(v)
\\\label{eq:t2} &&-\rho(x_2,x_4)\rho(x_1,x_3)(v)
+\rho(x_2,x_3)\rho(x_1,x_4)(v),\\
 \nonumber0&=&\rho(x_1,x_2)\rho(x_3,x_4)(v) - \rho(x_3,x_4)\rho(x_1,x_2)(v) -\rho([x_1,x_2,x_3],x_4)(v)
 \label{eq:t3} \\&&-   \nu(x_4)(v,\omega(x_1,x_2,x_3))-\rho(x_3,[x_1,x_2,x_4])(v)+ \nu(x_3)(v,\omega(x_1,x_2,x_4)).
  \end{eqnarray}

  \emptycomment{
  \begin{eqnarray}
\label{eq:t3} \nonumber&&\nu([x_1,x_2, x_3])(v_1,v_2)-\rho(x_2,x_3)\nu(x_1)(v_1,v_2)+\rho(x_1,x_3)\nu(x_2)(v_1,v_2)\\&&-\rho(x_1,x_2)\nu(x_3)(v_1,v_2)=0,\\
\label{eq:t4} \nonumber&&\rho(x_1,x_2)\nu(x_3)(v_1,v_2)-\nu(x_3)(\rho(x_1,x_2)(v_1),v_2)-\nu(x_3)(v_1,\rho(x_1,x_2)(v_2))\\&&-\nu([x_1,x_2,x_3])(v_1,v_2)=0,\\
\label{eq:t5} \nonumber&&\nu(x_1)(\rho(x_2,x_3)(v_2),v_1)+\rho(x_2,x_3)(\nu(x_1)(v_1,v_2))-\nu(x_3)(v_2,\rho(x_1,x_2)(v_1))\\
&&+\nu(x_2)(v_2,\rho(x_1,x_3)(v_1))=0,\\
\label{eq:t6}&&\nu(x_1)(\nu(x_2)(v_2,v_3),v_1)+\nu(x_2)(\nu(x_1)(v_1,v_2),v_3)+\nu(x_2)(v_2,\nu(x_1)(v_1,v_3))=0,\\
\label{eq:t7}&&\nu(x_2)(v_3,\nu(x_1)(v_1,v_2))+\nu(x_1)(\nu(x_2)(v_1,v_2),v_3)=0.
  \end{eqnarray}
}
\end{thm}
\pf The pair  $(\g\oplus V,[\cdot,\cdot,\cdot]_{(\rho,\nu,\omega)})$ defines a $3$-Lie algebra if and only if for all $e_1,\cdots,e_5\in\g\oplus V$, $$F_{e_1,e_2,e_3,e_4,e_5}=0.$$ By $F_{x_1,x_2,x_3,x_4,x_5}=0$, we have
\begin{eqnarray*}[x_1,x_2, [x_3,x_4,x_5]_{(\rho,\nu,\omega)}]_{(\rho,\nu,\omega)}-[[x_1,x_2, x_3]_{(\rho,\nu,\omega)},x_4,x_5]_{(\rho,\nu,\omega)}-[x_3,[x_1,x_2,x_4]_{(\rho,\nu,\omega)},x_5]_{(\rho,\nu,\omega)}\\
-[x_3,x_4,[x_1,x_2,x_5]_{(\rho,\nu,\omega)}]_{(\rho,\nu,\omega)}=0,\end{eqnarray*}
which gives Eq. \eqref{eq:t1}.

Similarly, $F_{x_1,v,x_2,x_3,x_4}=0$ gives Eq. \eqref{eq:t2}. $F_{x_1,x_2,v,x_3,x_4}=0$ gives Eq. \eqref{eq:t3}. $F_{x_1,x_2, v_1,v_2, x_3 }=0$ gives Eq. \eqref{eq:r3}. $F_{v_1,x_1, v_2,x_2, x_3 }=0$ gives Eq. \eqref{eq:r4}.  
$F_{v_1,x_1, v_2,v_3,x_2 }=0$ gives Eq. \eqref{eq:r6}. $F_{v_1,v_2,v_3, x_1,x_2 }=0$ gives Eq. \eqref{eq:r7}.

Conversely, if Eqs. \eqref{eq:r3}-\eqref{eq:r7} and Eqs. \eqref{eq:t1}-\eqref{eq:t3} hold, it is straightforward to see that for all $e_1,\cdots,e_5\in\g\oplus V$, $F_{e_1,e_2,e_3,e_4,e_5}=0.$ Thus, $(\g\oplus V,[\cdot,\cdot,\cdot]_{(\rho,\nu,\omega)})$ is a $3$-Lie algebra.  \qed

\begin{rmk}
For an abelian extension of $3$-Lie algebras, we have seen that $(V;\rho,\nu)$ is not a generalized representation, but the failure is controlled by $\omega$. This is totally different from the case of Lie algebras. For Lie algebras, abelian extension will give us a representation  and a $2$-cocycle. Only when we consider nonabelian extensions of Lie algebras (\cite{nonabelin cohomology of Lie}), the phenomenon that $\rho$ is not a representation will occur. But now, even for abelian extension of $3$-Lie algebras, this phenomenon has already occurred.
\end{rmk}

\begin{ex}{\rm
  Let $\g$ be the 3-dimensional 3-Lie algebra defined, with respect to a  basis $\{x_1,x_2,x_3\}$, by the skew-symmetric bracket $[x_1,x_2,x_3]=x_1$.

 Let $V$ be a 2-dimensional vector space and $\{v_1,v_2\}$ its basis. The following maps $\rho$, $\nu$ and $\omega$   define an abelian extension  of the 3-Lie algebra $\g$ on $V$ (according to Theorem \ref{thm:abelian ext}).

\begin{align*}
& \rho(x_1,x_2)(v_1)=r_1v_1+\frac{s_2 r_1}{s_1} v_2,  \; \;
\rho(x_1,x_2)(v_2)=\frac{s_1 r_2}{s_2} v_1+r_2v_2 ,   \; \\
 & \rho(x_1,x_3)(v_1)=\frac{s_3 r_1}{s_1} v_1+\frac{s_2 s_3 r_1}{s_1^2} v_2, \;\;
  \rho(x_1,x_3)(v_2)=\frac{s_3 r_2}{s_2} v_1+\frac{ s_3 r_2}{s_1} v_2, \; \\
& \rho(x_2,x_3)(v_1)=\frac{s_2 r_1 r_3}{s_1 r_2} v_1+\frac{s_2^2 r_1 r_3}{s_1^2r_2} v_2,    \; \;
  \rho(x_2,x_3)(v_2)=r_3 v_1+\frac{s_2 r_3}{s_1} v_2 ,
 \\
& \nu(x_1)(v_1,v_2)=0,  \;
\nu(x_2)(v_1,v_2)=s_1 v_1+s_2 v_2,  \;
\nu(x_3)(v_1,v_2)=s_3 v_1+\frac{s_2 s_3}{s_1} v_2,
\\
&\omega(x_1,x_2,x_3)=\frac{r_2}{s_2} v_1 -\frac{r_1}{s_1} v_2,
\end{align*}
where $r_i$ and $s_j$ are parameters in $\K$.

Notice that the map $\rho$ is not a representation in the usual sense (identities \eqref{eq:r1} and \eqref{eq:r2} are not satisfied). It becomes a generalized representation if and only if $r_1=r_2=0$.
}
\end{ex}

Now, we describe non-split abelian extensions using Maurer-Cartan elements.
The set $MC(L)$ of {\bf Maurer-Cartan elements} of a DGLA $(L,[\cdot,\cdot],\dM)$ is defined by
$$MC(L)\triangleq \{P\in L_1\mid\dM P+\frac{1}{2}[P,P]=0\}.$$

Let $(\g,[\cdot,\cdot,\cdot]_\g)$ be a 3-Lie algebra and $V$ a vector space. Let $\g\oplus V$ be the 3-Lie algebra direct sum of   $\g$ and $V$, where the bracket is defined by  $[x+u,y+v,z+w]=[x,y,z]_{\g}$. Then there is a DGLA   $(C(\g\oplus V,\g\oplus V), \nrn{\cdot,\cdot,\cdot},\delta)$, where $C(\g\oplus V,\g\oplus V)=\oplus_{p\geq 0}C^p(\g\oplus V,\g\oplus V)$ and $\delta$ is the coboundary operator for the 3-Lie algebra $\g\oplus V$ with the coefficients in the adjoint representation. 
  It is not difficult  to see that $(C_>(\g\oplus V, V),\nrn{\cdot,\cdot,\cdot},\delta)$ is a sub-DGLA of $(C(\g\oplus V,\g\oplus V), \nrn{\cdot,\cdot,\cdot},\delta)$,
  where $C_>(\g\oplus V, V)=\oplus_{p\geq0}C_>^p(\g\oplus V, V)$ is defined by \eqref{eq:C}.
\begin{pro}\label{pro:JMC}
The following two statements are equivalent:
\begin{itemize}
  \item[\rm(a)] $(\g\oplus V,[\cdot,\cdot,\cdot]_{(\rho,\nu,\omega)})$ is a $3$-Lie algebra, which is a non-split abelian extension of $\g$ by $V$;

  \item[\rm(b)] $\bar{\rho}+\bar{\nu}+\omega$ is a Maurer-Cartan element of the DGLA $(C_>(\g\oplus V, V),\nrn{\cdot,\cdot,\cdot},\delta)$.
\end{itemize}

\end{pro}
\pf By Lemma \ref{lem:cs}, $(\g\oplus V,[\cdot,\cdot,\cdot]_{(\rho,\nu,\omega)})$ is a 3-Lie algebra if and only if
$$\nrn{\pi+\bar{\rho}+\bar{\nu}+\omega,\pi+\bar{\rho}+\bar{\nu}+\omega}=0,$$
which can be rewritten as
 $$\nrn{\pi,\bar{\rho}+\bar{\nu}+\omega}+\half\nrn{\bar{\rho}+\bar{\nu}+\omega,\bar{\rho}+\bar{\nu}+\omega}=0.$$
 Since $\nrn{\pi,\bar{\rho}+\bar{\nu}+\omega}=\delta(\bar{\rho}+\bar{\nu}+\omega)$, we have
 $$\delta(\bar{\rho}+\bar{\nu}+\omega)+\half\nrn{\bar{\rho}+\bar{\nu}+\omega,\bar{\rho}+\bar{\nu}+\omega}=0,$$
which means that $\bar{\rho}+\bar{\nu}+\omega$ is a Maurer-Cartan element. \qed



\begin{thebibliography}{999}


\bibitem{AMS11} J. Arnlind, A. Makhlouf and S.  Silvestrov,  Construction of $n$-Lie algebras and $n$-ary Hom-Nambu-Lie algebras. \emph{J. Math. Phys. }\textbf{52} (2011), no. 12, 123502.



\bibitem{BL0}
J. Bagger and N. Lambert,  Gauge symmetry and supersymmetry of multiple M2-branes gauge theories. \emph{ Phys. Rev. D} \textbf{77} (2008), 065008.
\bibitem{BL3}
J. Bagger and N. Lambert,  Three-algebras and N=6 Chern-Simons gauge theories. \emph{ Phys. Rev. D} \textbf{79} (2009), no. 2, 025002, 8 pp.

\bibitem{Realization}
R. Bai, C. Bai and J. Wang, Realizations of 3-Lie algebras. \emph{J. Math. Phys.}, 51 (2010), 063505.

\bibitem{Baiclassification} R. Bai,  G. Song and Y. Zhang,   On classification of $n$-Lie algebras. \emph{Front. Math. China} \textbf{6} (2011), no. 4, 581--606.



\bibitem{Basu} A. Basu and J. A. Harvey, The M2-M5 brane system and a generalized Nahm's equation. \emph{Nucl. Phys. B 713}, 136 (2005), 136--150.







  \bibitem{DT} Y. Daletskii and L. Takhtajan, Leibniz and Lie algebra structures for Nambu
algebra. \emph{Lett. Math. Phys.} \textbf{39} (1997), 127--141.


 \bibitem{review}
J. A. de Azc$\rm\acute{a}$rraga and J. M. Izquierdo, $n$-ary algebras: a review with applications,
\emph{ J. Phys. A: Math. Theor.} \textbf{43} (2010), 293001.


\bibitem{cohomology}
J. A. de Azc$\rm\acute{a}$rraga and J. M. Izquierdo,
Cohomology of Filippov algebras and an
analogue of Whitehead's lemma.  	\emph{J. Phys. Conf. Ser.} \textbf{175}: 012001, (2009).


\bibitem{deformation}
J. Figueroa-O$'$Farrill, Deformations of 3-algebras. \emph{J. Math. Phys.} \textbf{50} (2009), no. 11, 113514, 27 pp.

\bibitem{Filippov}  V. T. Filippov, $n$-Lie algebras.  {\it Sib. Mat. Zh.} \textbf{26} (1985) 126--140.

\bibitem{nonabelin cohomology of Lie}
Y. Frégier, Non-abelian cohomology of extensions of Lie algebras as Deligne groupoid. \emph{J. Algebra} \textbf{398} (2014) 243--257.

\bibitem{Gautheron}
P. Gautheron, Some remarks concerning Nambu mechanics. \emph{Lett. Math. Phys.} \textbf{37} (1996) 103--116.


\bibitem{BL2}
J. Gomis, D. Rodríguez-Gómez, M. Van Raamsdonk and H. Verlinde,  Supersymmetric Yang-Mills theory from Lorentzian three-algebras. \emph{J. High Energy Phys.}, no. 8, 094, 18 pp  (2008).

\bibitem{HHM} P. Ho, R. Hou and Y. Matsuo, Lie $3$-algebra and multiple
$M_2$-branes.  \emph{J. High Energy Phys.},  no. 6, 020, 30 pp  (2008).

\bibitem{Kasymov}
Sh. M. Kasymov, On a theory of $n$-Lie algebras. (Russian) \emph{Algebra i Logika} \textbf{26} , no. 3 (1987) 277--297.

\bibitem{Makhlouf} A. Makhlouf, On Deformations of $n$-Lie Algebras, Chapter 4 in
 Non Associative \& Non Commutative Algebra
and Operator Theory, C.T. Gueye, M.S. Molina (eds.), Springer Proceedings in Mathematics \& Statistics \textbf{160}, (2016).

\bibitem{N} Y. Nambu, Generalized Hamiltonian dynamics. {\it Phys. Rev. D} \textbf{7} (1973)
                2405--2412.

 \bibitem{NR}
A. Nijenhuis and R. Richardson, Cohomology and Deformations in Graded Lie Algebras. \emph{Bull. Amer. Math. Soc.} \textbf{72} (1966) 1--29.


\bibitem{P} G. Papadopoulos, M2-branes, $3$-Lie algebras and
Plucker relations.  \emph{J. High Energy Phys.}  (2008),  no. 5, 054, 9 pp.

 \bibitem{NR bracket of n-Lie}
M. Rotkiewicz, Cohomology ring of $n$-Lie algebras. \emph{Extracta Math.} \textbf{20} no. 3, (2005) 219--232.


\bibitem{T} L. Takhtajan, On foundation of the generalized Nambu mechanics.
{\it Comm. Math. Phys.} \textbf{160} (1994) 295--315.
\bibitem{Tcohomology}
L. Takhtajan, A higher order analog of Chevalley-Eilenberg complex and deformation
theory of $n$-algebras. \emph{St. Petersburg Math. J.} \textbf{6} (1995) 429--438.



\end{thebibliography}
\end{document}